\documentclass[a4paper,oneside,12pt]{article}
\usepackage{amsmath, amsthm, amsfonts, amssymb}
\input xy   
\xyoption{all}
\usepackage{xcolor}

\numberwithin{equation}{section}
\swapnumbers

\theoremstyle{plain}
  \begingroup
        
        \newtheorem{theorem}[equation]{Theorem}
        \newtheorem{corollary}[equation]{Corollary}
        \newtheorem{proposition}[equation]{Proposition}
        \newtheorem{remark}[equation]{Remark}
        \newtheorem{lemma}[equation]{Lemma}
  \endgroup

\theoremstyle{definition}
  \begingroup
	\newtheorem{definition}[equation]{Definition}
        \newtheorem{notation}[equation]{Notation}

  \endgroup

\newcommand{\mr}[1]{\buildrel {#1} \over \longrightarrow}

\newcommand{\Mr}[1]{\buildrel {#1} \over \Longrightarrow}

\newcommand{\cc}{\mathcal}

\begin{document}

\title{A construction of 2-filtered bicolimits of categories}

\author{Eduardo J. Dubuc and Ross Street}
%{\scshape Eduardo J. Dubuc and Ross Street}

\date{\vspace{-5ex}}

\maketitle

\begin{center}
\`A la m\'emoire de Charles Ehresmann
\end{center}

%\begin{abstract}
\leftskip=1.1cm \rightskip=1.1cm  \small
\baselineskip=0.49cm

\noindent {\bf R\'esum\'e.} Nous d\'efinissons la notion de
\emph{2-cat\'egorie 2-filtrante} et donnons 
une construction explicite du bicolimite d'un 2-foncteur \`a valeurs
dans les cat\'egories.  
Une cat\'egorie consid\'er\'ee une  \mbox{2-cat\'egorie} triviale est
 2-filtrante si et seulement si c'est une cat\'egorie \mbox{filtrante,}
et notre construction conduit \`a une cat\'egorie \'equivalente
 \`a la cat\'egorie qui s'obtient par la
 construction usuelle de colimites filtrantes de
cat\'egories.
Pour cette construction des axiomes plus faibles suffisent, et nous
appelons la notion 
\mbox{correspondante} \mbox{\emph{2-cat\'egorie
pr\'e 2-filtrante}.}
L'ensemble complet des axiomes est n\'ecessaire 
pour montrer que les bicolimites 2-filtrantes 
ont les \mbox{propri\'et\'es} correspondantes aux propri\'et\'es
essentielles des colimites \mbox{filtrantes.}  

\leftskip=0pt \rightskip=0pt \normalsize
\baselineskip=0.49cm
%\end{abstract}

\vspace{2ex}

\noindent{\bf Introduction.}

\vspace{1ex}

We define the notion of \emph{2-filtered 2-category} and give an explicit
construction of the bicolimit of a category valued 2-functor. A category
considered as a trivial 2-category is 2-filtered if and only if it is a
filtered category, and our construction yields a category equivalent to
the category resulting from the usual construction of filtered
colimits of categories.  Weaker axioms 
suffice for this construction, and we call the corresponding notion
\emph{pre 2-filtered 2-category}. 
The full set of \mbox{axioms} is necessary to prove that 2-filtered
bicolimits have the properties
corresponding to the essential properties of filtered bicolimits.

In \cite{K} Kennison already considers filterness conditions on a
\mbox{2-category} under the name of \emph{bifiltered 2-category}.
It is easy to check that a bifiltered 2-category is 2-filtered, so our
results apply to bifiltered 2-categories.
Actually  Kennison's notion is equivalent to ours, but the other direction of
this equivalence is not entirely trivial.

\vspace{1ex}

\noindent {\bf Acknowledgement.} This paper gained its impetus from initial 
        collaboration in Montr\'eal
        amongst Andr\'e Joyal and the authors. We are very grateful for 
        Andr\'e's
        input. As authors of this account of that work and some further
        developments, we remain responsible for any imperfections.

%\vspace{8ex}

\section*{\normalsize 1. Pre 2-Filtered 2-Categories and the
  construction LL} 
\label{recall}
\addtocounter{section}{1}

%\section*{\normalsize 2. Weil prolongations of Banach manifolds.}
%\addtocounter{section}{1} \label{weilprolongations} 

\begin{definition} \label{pre_2}
A 2-category $\cc{A}$ is defined to be \emph{pre-2-filtered} when it
satisfies the following two axioms:

\vspace{1ex}

\noindent F1.
$$
Given \hspace{3ex}
\xymatrix@ur@R=2.5ex@C=2.5ex
         {
            E \ar[r]^f   \ar[d]_g 
          & A
          \\
            B   
         }
 \hspace{4ex} there \; exists \; invertible \hspace{2ex}
\xymatrix@R=2.5ex@C=2.5ex@ur
         {
            E \ar[r]^f  \ar@{}[dr]|{\gamma \Downarrow}  \ar[d]_g 
          & A \ar[d]^u
          \\
            B  \ar[r]_v 
          & C 
         }\;.
$$
\noindent F2.
$$
 Given \; any \; 2\,{\text-}\,cells \hspace{3ex}
\xymatrix@ur@R=2.5ex@C=2.5ex
         {
            E \ar[r]^f  \ar@{}[dr]|{\gamma_1 \Downarrow}  \ar[d]_g 
          & A \ar[d]^{u_1}
          \\
            B  \ar[r]_{v_1} 
          & C_1 
         }
 ,  \;
\xymatrix@ur@R=2.5ex@C=2.5ex
         {
            E \ar[r]^f  \ar@{}[dr]|{\gamma_2 \Downarrow}  \ar[d]_g 
          & A \ar[d]^{u_2}
          \\
            B  \ar[r]_{v_2} 
          & C_2 
         }
$$    
$$
there \; exists  
\xymatrix@R=2.5ex@C=2.5ex@ur
         {
          & C_1 \ar[d]^{w_1}
          \\
            C_2  \ar[r]_{w_2} 
          & C 
         }
\;\;\;\; with \; invertible \; 2\,{\text-}\,cells \;\; \alpha, \; \beta\; 
\hspace{1ex} 
such \; that 
$$
\begin{equation} \label{ll}
\xymatrix@ur@R=2.5ex@C=2.5ex
         {
            E \ar[r]^f  \ar@{}[dr]|{\gamma_1 \Downarrow}  \ar[d]_g 
          & A \ar[d]^{u_1}
          \\
            B \ar[r]_{\!\!v_1}  \ar@{}[dr]|{\alpha \Downarrow}  \ar[d]_{v_2} 
          & C_1 \ar[d]^{w_1}
          \\
            C_2 \ar[r]_{w_2}
          & C
         }
\hspace{4ex} = \hspace{4ex}
\xymatrix@ur@R=2.5ex@C=2.5ex
         {
            E \ar[r]^f  \ar@{}[dr]|{\gamma_2  \Downarrow}  \ar[d]_g 
          & A \ar[r]^{u_1}  \ar@{}[dr]|{\beta \Downarrow}
	    \ar[d]^{\!\!\!\!u_2} 
          & C_1 \ar[d]^{w_1}
          \\
            B  \ar[r]_{v_2} 
          & C_2 \ar[r]_{w_2}
          & C
         }
\;\;\;\;.
\end{equation}

\end{definition}
\begin{notation}[the LL equation]
Given two pairs of 2-cells $(\gamma_1, \; \alpha)$,
$(\gamma_2, \; \beta)$ as in F2,  we shall call the equation \ref{ll} the
\mbox{\emph{equation LL}}. 
\end{notation}

%\pagebreak

Axioms F1 and F2 can be weakened if we add a third axiom: 

\vspace{1ex}

\noindent WF1. Same axiom F1 but not requiring $\gamma$ to be invertible.

\vspace{1ex}

\noindent WF2. Same axiom F2 but only for invertible $\gamma_1$
and $\gamma_2$.

\vspace{1ex}

\noindent WF3.
$$
Given \hspace{1ex}
\xymatrix@R=2.5ex@C=2.5ex@ur
         {
            E \ar[r]^f  \ar@{}[dr]|{\gamma \Downarrow}  \ar[d]_g 
          & A \ar[d]^u
          \\
            B  \ar[r]_v 
          & C 
         }
\hspace{2ex} there \; exists \hspace{1ex}
invertible \; 2\,{\text-}\,cells \;\; \alpha, \; \beta\;
such \; that 
$$
$$
\xymatrix@ur@R=2.5ex@C=2.5ex
         {
            E \ar[r]^f  \ar@{}[dr]|{\gamma \Downarrow}  \ar[d]_g 
          & A \ar[d]^{u}
          \\
            B \ar[r]_{\!\!v}  \ar@{}[dr]|{\alpha \Downarrow}  \ar[d]_{s_2} 
          & C \ar[d]^{h_2}
          \\
            C_2 \ar[r]_{w_2}
          & D_2
         }
\hspace{2ex} and \hspace{2ex}
\xymatrix@ur@R=2.5ex@C=2.5ex
         {
            E \ar[r]^f  \ar@{}[dr]|{\gamma  \Downarrow}  \ar[d]_g 
          & A \ar[r]^{s_1}  \ar@{}[dr]|{\beta \Downarrow}
	    \ar[d]^{\!\!\!\!u} 
          & C_1 \ar[d]^{w_1}
          \\
            B  \ar[r]_{v} 
          & C \ar[r]_{h_1}
          & D_1
         }
\hspace{4ex} are \; invertible.
$$ 

\vspace{1ex}

We leave to the reader the proof of the following:
\begin{proposition}
The set of axioms F1, F2 is equivalent  to the set WF1 WF2 and WF3.
\qed \end{proposition}
When $\cc{A}$ is a trivial 2-category (the only 2-cells are the
identities), our axiom F1 corresponds to axiom PS1 in the definition of
\mbox{pseudofiltered} category (cf \cite{G1} Expos\'{e} I), while axiom PS2 may
not hold. Thus, 
a category which is pre 2-filtered as a trivial 2-category may not be
\mbox{pseudofiltered.} Notice
that our axioms F2 and WF3 are vacuous in this case.
 
\vspace{1ex}

%\section*{\normalsize 2. The construction LL}
%\addtocounter{section}{1} \label{constructionLL} 

{\bf Construction LL}

Let $\cc{A}$ be a pre-2-filtered 2-category and $F: \cc{A} \to
\cc{C}at$ a category valued \mbox{2-functor.} 
We shall now construct a category
which is to be the bicolimit (in the sense made precise in
theorem \ref{bicolimit}) \mbox{of $F$.}
This construction generalizes 
Grothendieck's construction of the   
%$\stackrel{\txt{\normalsize{Lim F}}}
%{\stackrel{\!\!\longrightarrow}{\txt{\scriptsize {$\cc{A}$}}}}$
%
\mbox{category
$\stackrel{\stackrel{\txt{\normalsize{Lim F}}}{\!\!\!\!\longrightarrow}}
{\txt{\scriptsize {$\!\!\!\! \cc{A}$}}}$ for a filtered
category $\cc{A}$} (cf \cite{G2}, Expos\'{e} VI).   

\vspace{1ex}

\begin{definition}[Quasicategory $\cc{L}(F)$] \label{constructionLL} $\\$
$\hspace{2.2ex}$
i) An \emph{object} is a pair $(x,\, A)$ with $x \in FA$. 

\vspace{1ex}

ii) A \emph{premorphism}  $(x,\, A) \to (y,\, B)$ between two objects
 is a triple $(u,\, \xi, \, v)$, where 
 and  $A \mr{u} C$, 
 $B \mr{v} C$ and  
\mbox{$\xi: F(u)(x) \to F(v)(y)$ in $FC$.}

\vspace{1ex}

iii) A \emph{homotopy} between two premorphisms is a quadruple 
$(w_1, \; w_2, \; \alpha, \; \beta): ({u_1},\, {\xi_1}, \, {v_1})
\Rightarrow ({u_2},\, {\xi_2}, \, {v_2})$, where $C_1 \mr{w_1} C$, 
 $C_2 \mr{w_2} C$ and  
 \mbox{$\alpha : {w_1}{v_1} \mr{\cong} {w_2}{v_2}$,} 
 \mbox{$\beta : {w_1}{u_1} \mr{\cong} {w_2}{u_2}$} are invertible
 \mbox{2-cells} such that
$$
\xymatrix@C=7ex
          {
           F(w_1)F(u_1)(x) = F(w_1u_1)(x) \ar[r]^{F(\beta)x}
                                          \ar@<-6ex>[d]^{F(w_1)(\xi_1)}
          &
           F(w_2u_2)(x) = F(w_2)F(u_2)(x) \ar@<6ex>[d]^{F(w_2)(\xi_2)}
          \\
           F(w_1)F(v_1)(y) = F(w_1v_1)(y) \ar[r]^{F(\alpha)y}
          &
           F(w_2v_2)(y) = F(w_2)F(v_2)(y)
          }
$$
commutes in $FC$. 
\end{definition}

\vspace{1ex}

We shall formally introduce now an \emph{abuse of notation}
\begin{notation} \label{abuse} $\\$
$\hspace{2ex}$ 
i) We omit the letter $F$ in denoting the action of $F$ on
its arguments. Thus 
$
\xymatrix
         {
          A \ar@<1.5ex>[r]^u 
          \ar@{}@<-1ex>[r]^{\alpha \, \Downarrow} \ar@<-1.2ex>[r]_v & B
         }
$
indicates a 2-cell in $\cc{A}$ as well as the corresponding
natural transformation 
$
\xymatrix@C=8ex
         {
          FA \ar@<1.5ex>[r]^{F(u)} 
      \ar@{}@<-1ex>[r]^{F(\alpha) \, \Downarrow} \ar@<-1.2ex>[r]_{F(v)} & FB
         }
$
in $\cc{C}at$. In this way, the above commutative square becomes
$$
\xymatrix
          {
           w_1u_1x \ar[r]^{\beta x} \ar[d]^{w_1\xi_1}
          &
           w_2u_2x  \ar[d]^{w_2\xi_2}
          \\
           w_1v_1y  \ar[r]^{\alpha y}
          &
           w_2v_2y 
          }
$$
\hspace{2ex} ii) We write $F \mr{x} A$
in $(\cc{C}at^{\cc{A}})^{op}$ 
for the natural transformation \mbox{$\cc{A}[A, \, -] \to F$} defined by 
\mbox{$x_{C}(A \mr{u} C) = F(u)(x) \in FC$.} 

Notice that the notation in i) $F(u)(x) = ux$ is consistent with this
since juxtaposition denotes composition. Also, in the same vein, given
an object $x \in FA$ and a functor $FA \mr{h} \cc{X}$ into any other
category, the composite $F \mr{x} A \mr{h} \cc{X}$ makes sense and we
have $h(F(x)) = hx$.

\end{notation}

In this notation then, a premorphism  
$(x,\, A) \to (y,\, B)$is a triple $(u,\, \xi, \, v)$, where
$
\xymatrix@R=2.5ex@C=2.5ex@ur
         {
            F \ar[r]^x  \ar@{}[dr]|{\xi \Downarrow}  \ar[d]_y 
          & A \ar[d]^u
          \\
            B  \ar[r]_v 
          & C 
         }$;
that is  \mbox{$\xi : ux \mr{} vy$ in FC.} 
A homotopy between two premorphisms 
$\xymatrix@ur@R=2.5ex@C=2.5ex
         {
            F \ar[r]^x  \ar@{}[dr]|{\xi_1 \Downarrow}  \ar[d]_y 
          & A \ar[d]^{u_1}
          \\
            B  \ar[r]_{v_1} 
          & C_1
         }, \;
\xymatrix@ur@R=2.5ex@C=2.5ex
         {
            F \ar[r]^x  \ar@{}[dr]|{\xi_2 \Downarrow}  \ar[d]_y 
          & A \ar[d]^{{u_2}}
          \\
            B  \ar[r]_{{v_2}} 
          & C_2
         }$
is a pair of invertible 2-cells
$
\xymatrix@ur@R=2.5ex@C=2.5ex
         {
          B \ar[r]^{v_1}  \ar@{}[dr]|{\alpha \Downarrow}
	                                                    \ar[d]_{{v_2}}    
          & C_1 \ar[d]^{w_1}
          \\
            C_2  \ar[r]_{{w_2}} 
          & C
         } , \;
\xymatrix@ur@R=2.5ex@C=2.5ex
         {
          A \ar[r]^{u_1}  \ar@{}[dr]|{\beta \Downarrow}
	                                                    \ar[d]_{{u_2}}    
          & C_1 \ar[d]^{w_1}
          \\
            C_2  \ar[r]_{{w_2}} 
          & C
         } 
$
satisfying \mbox{the LL equation:}
\mbox{\hspace{6ex}     
$\xymatrix@ur@R=2.5ex@C=2.5ex
         {
            F \ar[r]^x  \ar@{}[dr]|{\xi_1 \Downarrow}  \ar[d]_y 
          & A \ar[d]^{u_1}
          \\
            B \ar[r]_{\!\!{v_1}}  \ar@{}[dr]|{\alpha \Downarrow}
	                                                     \ar[d]_{{v_2}}  
          & C_1 \ar[d]^{w_1}
          \\
            C_2 \ar[r]_{{w_2}}
          & C
         }
\hspace{1ex} = \hspace{1ex}
\xymatrix@ur@R=2.5ex@C=2.5ex
         {
            F \ar[r]^x  \ar@{}[dr]|{\xi_2 \Downarrow}  \ar[d]_y 
          & A \ar[r]^{u_1}  \ar@{}[dr]|{\beta \Downarrow}
	                                             \ar[d]^{\!\!\!\!{u_2}} 
          & C_1 \ar[d]^{w_1}
          \\
            B  \ar[r]_{{v_2}} 
          & C_2 \ar[r]_{{w_2}}
          & C
         }.
$} 

\vspace{1ex}

We shall simply write 
$(\alpha,\, \beta): \xi_1 \Rightarrow \xi_2$ for all the data involved
in a homotopy.

At this point it is convenient to introduce the following notation:

%\pagebreak

\begin{notation} [LL-composition of 2-cells] \label{composite}
Given three 2-cells \mbox{$\alpha$, $\beta$, and $\gamma$} that fit into a
diagram as it follows, we write:
$$
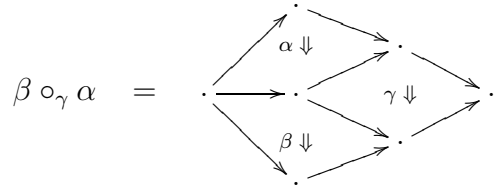
\xymatrix@R=1.5ex@C=4ex
           {
            \\
            \\
              \beta \circ_\gamma  \alpha  \hspace{2ex} =  \hspace{2ex}
            }
\xymatrix@R=1ex@C=4ex
         {
          &
           . \ar[rd]
          \\
          &
           {\txt{\scriptsize{$\alpha \Downarrow$}}}
          &
           . \ar[rd]
          \\
           . \ar[r] \ar[rdd] \ar[ruu]
          &
           . \ar[ru] \ar[rd] 
          &
           {\txt{\scriptsize{$\gamma \Downarrow$}}}
          &
           .
          \\
          &
           {\txt{\scriptsize{$\beta \Downarrow$}}}
          &
           . \ar[ru]
          \\
          &
           . \ar[ru]
         }    
$$
Thus, $\beta \circ_\gamma  \alpha$ is our notation for the 2-cell between the
top and the bottom composites of arrows. It should be thought of as the
``composite of \mbox{$\beta$ with $\alpha$ over $\gamma$''.}
\end{notation}

\noindent Homotopies compose: \mbox{Consider a third premorphism 
$\xymatrix@ur@R=2.5ex@C=2.5ex
         {
            F \ar[r]^x  \ar@{}[dr]|{\xi_3 \Downarrow}  \ar[d]_y 
          & A \ar[d]^{{u_3}}
          \\
            B  \ar[r]_{{v_3}} 
          & C_3
         }$}
\mbox
     {
      and a homotopy
$
\xymatrix@ur@R=2.5ex@C=2.5ex
         { 
          B \ar[r]^{v_2}  \ar@{}[dr]|{\alpha' \Downarrow}
	                                                       \ar[d]_{{v_3}}             & C_2 \ar[d]^{w'_2}
          \\
            C_3  \ar[r]_{{w'_3}} 
          & C'
         } , \;
\xymatrix@ur@R=2.5ex@C=2.5ex
         {
          A \ar[r]^{u_2}  \ar@{}[dr]|{\beta' \Downarrow}
	                                                     \ar[d]_{{u_3}}    
          & C_2 \ar[d]^{w'_2}
          \\
            C_3  \ar[r]_{{w'_3}} 
          & C'
         }
$, 
$(\alpha',\, \beta'): \xi_2 \Rightarrow \xi_3$.  Use 
    } 
axiom F1 to obtain an invertible 2-cell
$\xymatrix@ur@R=2.5ex@C=2.5ex
         {
          C_2 \ar[r]^{w_2}  \ar@{}[dr]|{\gamma \Downarrow}
                                                        \ar[d]_{{w'_2}}    
          & C \ar[d]^{h}
          \\
            C'  \ar[r]_{{h'}} 
          & H
         }$.
 This determines a pair of 2-cells:
$$
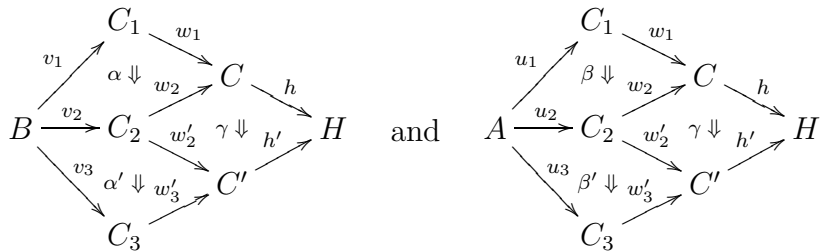
\hspace{2ex}
\xymatrix@R=1ex@C=4ex
         {
          &
           C_1 \ar[rd]^{w_1}
          \\
          &
           {\txt{\scriptsize{$\alpha \Downarrow$}}}
          &
           C \ar[rd]^{h}
          \\
           B \ar[r]^{v_2} \ar[rdd]^{\!\! v_3} \ar[ruu]^{v_1}
          &
           C_2 \ar[ru]^{w_2} \ar[rd]^{\!\! w'_2} 
          &
           {\txt{\scriptsize{$\gamma \Downarrow$}}}
          &
           H 
          \\
          &
           {\txt{\scriptsize{$\alpha' \Downarrow$}}}
          &
           C' \ar[ru]^{h'}
          \\
          &
           C_3 \ar[ru]^{w'_3}
         }
\hspace{2ex}
\xymatrix@R=6ex
           {
            \\
            \txt{and}
            }
\hspace{2ex}
\xymatrix@R=1ex@C=4ex
         {
          &
           C_1 \ar[rd]^{w_1}
          \\
          &
           {\txt{\scriptsize{$\beta \Downarrow$}}}
          &
           C \ar[rd]^{h}
          \\
           A \ar[r]^{u_2} \ar[rdd]^{\!\! u_3} \ar[ruu]^{u_1}
          &
           C_2 \ar[ru]^{w_2} \ar[rd]^{\!\! w'_2} 
          &
           {\txt{\scriptsize{$\gamma \Downarrow$}}}
          &
           H 
          \\
          &
           {\txt{\scriptsize{$\beta' \Downarrow$}}}
          &
           C' \ar[ru]^{h'}
          \\
          &
           C_3 \ar[ru]^{w'_3}
         }
$$

%\pagebreak

\noindent which defines a homotopy $\xi_1 \Rightarrow \xi_3$.
The corresponding LL equation follows easily from the LL equations for
$(\alpha, \;\beta)$ and $(\alpha', \; \beta')$.

\vspace{1ex}

Using notation \ref{composite}, we have:
\begin{proposition}[vertical composition of homotopies] \label{v.c.o.h}
Given
$(\alpha,\, \beta): \xi_1 \Rightarrow \xi_2$ and $(\alpha',\, \beta'): \xi_2
\Rightarrow \xi_3$,
 there exists an appropriate $\gamma$ and a homotopy 
$(\alpha' \circ_\gamma \alpha, \, \beta' \circ_\gamma \beta): \xi_1
\Rightarrow \xi_3$. 
\qed \end{proposition} 
 
\vspace{1ex}

Homotopies are generated by composition out of two basic ones. The
proof of the following is immediate:

\begin{proposition} \label{basich}
Every pair of premorphisms $\xi_1$, $\xi_2$ and pair of 2-cells
$\alpha$, $\beta$ that fit as follows determine two basic
\mbox{homotopies:}  
$$
(\alpha, \; id_1) : 
\xymatrix@ur@R=2.5ex@C=2.5ex
         {
            F \ar[r]^x  \ar@{}[dr]|{\xi_1 \Downarrow}  \ar[d]_y 
          & A \ar[d]^{u_1}
          \\
            B  \ar[r]_{v_1} 
          & C_1
         }
\hspace{4ex} \Rightarrow \hspace{4ex}
\xymatrix@ur@R=2.5ex@C=2.5ex
         {
            F \ar[r]^x  \ar@{}[dr]|{\xi_1 \Downarrow}  \ar[d]_y 
          & A \ar[d]^{u_1}
          \\
            B \ar[r]_{\!\!{v_1}}  \ar@{}[dr]|{\alpha \Downarrow}
	                                                     \ar[d]_{{v_2}}  
          & C_1 \ar[d]^{w_1}
          \\
            C_2 \ar[r]_{{w_2}}
          & C
         }
$$
$$
(id_2, \; \beta):
\xymatrix@ur@R=2.5ex@C=2.5ex
         {
            F \ar[r]^x  \ar@{}[dr]|{\xi_2 \Downarrow}  \ar[d]_y 
          & A \ar[r]^{u_1}  \ar@{}[dr]|{\beta \Downarrow}
	                                             \ar[d]^{\!\!\!\!{u_2}} 
          & C_1 \ar[d]^{w_1}
          \\
            B  \ar[r]_{{v_2}} 
          & C_2 \ar[r]_{{w_2}}
          & C
         }
\hspace{4ex} \Rightarrow \hspace{4ex}
\xymatrix@ur@R=2.5ex@C=2.5ex
         {
            F \ar[r]^x  \ar@{}[dr]|{\xi_2 \Downarrow}  \ar[d]_y 
          & A \ar[d]^{{u_2}}
          \\
            B  \ar[r]_{{v_2}} 
          & C_2
         }
$$
where $id_1$ and $id_2$ are the identity 2-cells corresponding to the
arrows $w_1 u_1$ and $w_2 v_2$ respectively. When the pair $(\alpha,
\; \beta)$ satisfies the LL equation, then the composite
(over the identity 2-cell of the identity arrow of $C_2$) of these
basic homotopies is defined, and it is equal to the homotopy determined by
$(\alpha, \; \beta)$.
\qed \end{proposition}              

\pagebreak

Premorphisms compose: Given two premorphisms 

$
\xymatrix{
     (x,\,A) \ar[r]^{\xi} %\ar@/_2pc/[rr]^{(s,\, \mu, \, t)}
        & (y,\,B) \ar[r]^{\zeta}
        & (z,\,C)
          }
$, 
$
\xymatrix@R=2.5ex@C=2.5ex@ur
         {
            F \ar[r]^x  \ar@{}[dr]|{\xi \Downarrow}  \ar[d]_y 
          & A \ar[d]^u
          \\
            B  \ar[r]_v 
          & S 
         }
$,
$
\xymatrix@R=2.5ex@C=2.5ex@ur
         {
            F \ar[r]^y  \ar@{}[dr]|{\zeta \Downarrow}  \ar[d]_z 
          & B \ar[d]^h
          \\
            C  \ar[r]_k 
          & T 
         }
$, use axiom F1 to obtain invertible
$\xymatrix@ur@R=2.5ex@C=2.5ex
         {
          B \ar[r]^{v}  \ar@{}[dr]|{\gamma \Downarrow}
                                                        \ar[d]_{{h}}    
          & S \ar[d]^{s}
          \\
            T  \ar[r]_{{t}} 
          & H
         }.$ 
According to \ref{composite} this \mbox{determines} a premorphism $\zeta
\circ_\gamma \xi$ between $(x,\,A)$ and $(z,\,C)$ that we take as a
composite of $\xi$ with $\zeta$. Thus: 
$$
\xymatrix@R=1.5ex@C=4ex
           {
            \\
            \\
              \zeta \circ_\gamma  \xi  \hspace{2ex} =  \hspace{2ex}
            }
\xymatrix@R=1ex@C=4ex
         {
          &
           A \ar[rd]^u
          \\
          &
           {\txt{\scriptsize{$\xi \Downarrow$}}}
          &
           S \ar[rd]^s
          \\
           F \ar[r]^y \ar[rdd]_z \ar[ruu]^x
          &
           B \ar[ru]^v \ar[rd]^h 
          &
           {\txt{\scriptsize{$\gamma \Downarrow$}}}
          &
           H
          \\
          &
           {\txt{\scriptsize{$\zeta \Downarrow$}}}
          &
           T \ar[ru]_t
          \\
          &
           C \ar[ru]_k
         }    
$$
We have:  
\begin{proposition}[horizontal composition of premorphisms]
Given $\xi: (x,\,A) \to (y,\,B)$ and $\zeta: (y,\,B) \to (z,\,C)$, 
 there exists an appropriate $\gamma$ and
$\zeta \circ_\gamma \xi: (x,\,A) \to
(z,\,C)$. 
\qed \end{proposition}
Homotopies also compose horizontally:
\begin{proposition}[\mbox{horizontal composition of homotopies}]
\label{h.c.o.h}
Consider composable premorphisms and homotopies as follows: 
\xymatrix@C=8ex{
          {(x,\, A)} 
          \ar@<2ex>[r]^{\xi_1} \ar@{}@<-1.45ex>[r]^{(\alpha, \beta)
          \Downarrow}  
          \ar@<-2ex>[r]_{\xi_2}
        & {(y,\, B)}
          \ar@<2ex>[r]^{\zeta_1} \ar@{}@<-1.45ex>[r]^{(\varepsilon, \delta)
          \Downarrow}  
          \ar@<-2ex>[r]_{\zeta_2}
        & {(z,\, C)}
          }.
Then, given any two \mbox{composites} 
$\zeta_1 \circ_{\gamma_1}  \xi_1$ and $\zeta_2 \circ_{\gamma_2} \xi_2$,   
there exists a homotopy
$\zeta_1 \circ_{\gamma_1}  \xi_1  \Rightarrow \zeta_2 \circ_{\gamma_2} \xi_2$. 
\end{proposition}
\begin{proof}

\pagebreak

Consider the given homotopies and their LL equations:
$$
\xymatrix@ur@R=2.5ex@C=2.5ex
         {
            {F} \ar[r]^x  \ar@{}[dr]|{\xi_1 \Downarrow}  \ar[d]_y 
          & {A} \ar[d]^{u_1}
          \\
            {B} \ar[r]_{v_1}  \ar@{}[dr]|{\alpha \Downarrow} \ar[d]_{v_2}  
          & {S_1} \ar[d]^{s_1}
          \\
            {S_2} \ar[r]_{s_2}
          & {S}
         }
\hspace{4ex} = \hspace{4ex}
\xymatrix@ur@R=2.5ex@C=2.5ex
         {
            {F} \ar[r]^x  \ar@{}[dr]|{\xi_2 \Downarrow}  \ar[d]_y 
         & {A} \ar[r]^{u_1}  \ar@{}[dr]|{\beta \Downarrow}  \ar[d]^{\!\!\!u_2} 
          & {S_1} \ar[d]^{s_1}
          \\
            {B}  \ar[r]_{v_2} 
          & {S_2} \ar[r]_{s_2}
          & {S}
         }
$$
$$
\xymatrix@ur@R=2.5ex@C=2.5ex
         {
            {F} \ar[r]^y  \ar@{}[dr]|{\zeta_1 \Downarrow}  \ar[d]_z 
          & {B} \ar[d]^{h_1}
          \\
       {C} \ar[r]_{\!\!k_1}  \ar@{}[dr]|{\varepsilon \Downarrow} \ar[d]_{k_2}  
          & {T_1} \ar[d]^{t_1}
          \\
            {T_2} \ar[r]_{t_2}
          & {T}
         }
\hspace{4ex} = \hspace{4ex}
\xymatrix@ur@R=2.5ex@C=2.5ex
         {
            {F} \ar[r]^y  \ar@{}[dr]|{\zeta_2 \Downarrow}  \ar[d]_z 
        & {B} \ar[r]^{h_1}  \ar@{}[dr]|{\delta \Downarrow}  \ar[d]^{\!\!\!h_2} 
          & {T_1} \ar[d]^{t_1}
          \\
            {C}  \ar[r]_{k_2} 
          & {T_2} \ar[r]_{t_2}
          & {T}
         }
$$
and consider the horizontal composites of the premorphisms:
$$
\xymatrix@R=1ex@C=4ex
         {
          &
           A \ar[rd]^{u_1}
          \\
          &
           {\txt{\scriptsize{$\xi_1 \!\!\Downarrow$}}}
          &
           S_1 \ar[rd]^{r_1}
          \\
           F \ar[r]^y \ar[rdd]_z \ar[ruu]^x
          &
           B \ar[ru]^{v_1} \ar[rd]^{h_1} 
          &
           {\txt{\scriptsize{$\gamma_1 \!\!\Downarrow$}}}
          &
           H_1 
          \\
          &
           {\txt{\scriptsize{$\zeta_1 \!\!\Downarrow$}}}
          &
           T_1 \ar[ru]_{l_1}
          \\
          &
           C \ar[ru]_{k_1}
         }
\xymatrix@R=1.5ex
           {
            \\
            \\
            \\
            ,
            }
\hspace{2ex}
\xymatrix@R=1ex@C=4ex
         {
          &
           A \ar[rd]^{u_2}
          \\
          &
           {\txt{\scriptsize{$\xi_2 \!\!\Downarrow$}}}
          &
           S_2 \ar[rd]^{r_2}
          \\
           F \ar[r]^y \ar[rdd]_z \ar[ruu]^x
          &
           B \ar[ru]^{v_2} \ar[rd]^{h_2} 
          &
           {\txt{\scriptsize{$\gamma_2 \!\!\Downarrow$}}}
          &
           H_2
          \\
          &
           {\txt{\scriptsize{$\zeta_2 \!\!\Downarrow$}}}
          &
           T_2 \ar[ru]_{l_2}
          \\
          &
           C \ar[ru]_{k_2}
         }
$$
We shall produce a homotopy between these composites.

\vspace{1ex}

%\pagebreak

%First compute composits $\epsilon \circ_{\phi_1} \gamma_1$ and
%$\gamma_2 \circ_{\phi_2} \beta$. That is 

First use axiom F1 to obtain $\phi_1$, $\phi_2$ as follows:

$$
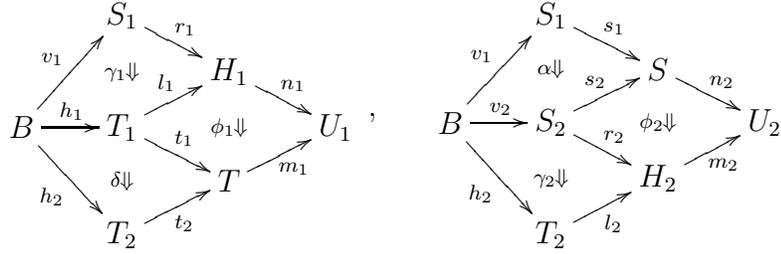
\hspace{2ex}
\xymatrix@R=1ex@C=4ex
         {
          &
           S_1 \ar[rd]^{r_1}
          \\
          &
           {\txt{\scriptsize{$\gamma_1 \!\!\Downarrow$}}}
          &
           H_1 \ar[rd]^{n_1}
          \\
           B \ar[r]^{h_1} \ar[rdd]_{h_2} \ar[ruu]^{v_1}
          &
           T_1 \ar[ru]^{l_1} \ar[rd]^{t_1} 
          &
           {\txt{\scriptsize{$\phi_1 \!\!\Downarrow$}}}
          &
           U_1 
          \\
          &
           {\txt{\scriptsize{$\delta \!\!\Downarrow$}}}
          &
           T \ar[ru]_{m_1}
          \\
          &
           T_2 \ar[ru]_{t_2}
         }
\xymatrix@R=6ex
           {
            \\
            , \;\;
            }
\hspace{2ex}
\xymatrix@R=1ex@C=4ex
         {
          &
           S_1 \ar[rd]^{s_1}
          \\
          &
           {\txt{\scriptsize{$\alpha \!\!\Downarrow$}}}
          &
           S \ar[rd]^{n_2}
          \\
           B \ar[r]^{v_2} \ar[rdd]_{h_2} \ar[ruu]^{v_1}
          &
           S_2 \ar[ru]^{s_2} \ar[rd]^{r_2} 
          &
           {\txt{\scriptsize{$\phi_2 \!\!\Downarrow$}}}
          &
           U_2
          \\
          &
           {\txt{\scriptsize{$\gamma_2 \!\!\Downarrow$}}}
          &
           H_2 \ar[ru]_{m_2}
          \\
          &
           T_2 \ar[ru]_{l_2}
         }
$$
Then use axiom F2 to obtain $\theta_1$, $\theta_2$ satisfying the LL
equation: 
$$
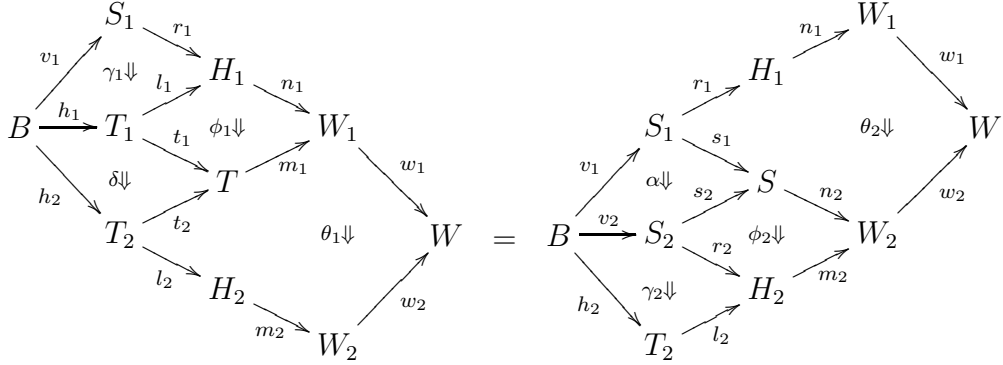
\xymatrix@R=1ex@C=4ex
         {
          &
           S_1 \ar[rd]^{r_1}
          \\
          &
           {\txt{\scriptsize{$\gamma_1 \!\!\Downarrow$}}}
          &
           H_1 \ar[rd]^{n_1}
          \\
           B \ar[r]^{h_1} \ar[rdd]_{h_2} \ar[ruu]^{v_1}
          &
           T_1 \ar[ru]^{l_1} \ar[rd]^{t_1} 
          &
           {\txt{\scriptsize{$\phi_1 \!\!\Downarrow$}}}
          &
           W_1 \ar[rdd]^{w_1} 
          \\
          &
           {\txt{\scriptsize{$\delta \!\!\Downarrow$}}}
          &
           T \ar[ru]_{m_1}
          \\
          &
           T_2 \ar[ru]_{t_2} \ar[rd]_{l_2}
          &
          &
           {\txt{\scriptsize{$\theta_1  \!\!\Downarrow$}}}
          &
           W
          \\
          &
          &
           H_2 \ar[rd]_{m_2}
          \\
          &
          &
          &
           W_2 \ar[ruu]_{w_2}
         }
\hspace{1ex}
\xymatrix@R=15ex
           {
            \\
            =
            }
\hspace{1ex}
\xymatrix@R=1ex@C=4ex
         {
          &
          &
          &
           W_1 \ar[rdd]^{w_1}
          \\
          &
          &
           H_1 \ar[ru]^{n_1} 
          \\
          &
           S_1 \ar[rd]^{s_1} \ar[ru]^{r_1}
          &
          & 
           {\txt{\scriptsize{$\theta_2 \!\!\Downarrow$}}}
          &
           W
          \\
          &
           {\txt{\scriptsize{$\alpha \!\!\Downarrow$}}}
          &
           S \ar[rd]^{n_2}
          \\
           B \ar[r]^{v_2} \ar[rdd]_{h_2} \ar[ruu]^{v_1}
          &
           S_2 \ar[ru]^{s_2} \ar[rd]^{r_2} 
          &
           {\txt{\scriptsize{$\phi_2 \!\!\Downarrow$}}}
          &
           W_2 \ar[ruu]_{w_2}
          \\
          &
           {\txt{\scriptsize{$\gamma_2 \!\!\Downarrow$}}}
          &
           H_2 \ar[ru]_{m_2}
          \\
          &
           T_2 \ar[ru]_{l_2}
         }
$$
The homotopy is given by the following pair of 2-cells:
$$
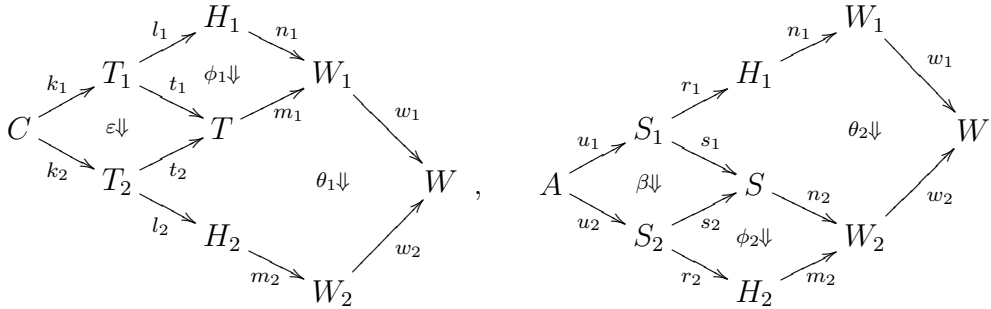
\xymatrix@R=1ex@C=4ex
         {
          &
          &
           H_1 \ar[rd]^{n_1}
          \\
          &
           T_1 \ar[ru]^{l_1} \ar[rd]^{t_1}                    
          &
           {\txt{\scriptsize{$\phi_1 \!\!\Downarrow$}}}
          &
           W_1 \ar[rdd]^{w_1}
          \\
           C \ar[ru]^{k_1} \ar[rd]_{k_2} 
          &
           {\txt{\scriptsize{$\varepsilon \!\!\Downarrow$}}} 
          &
           T \ar[ru]_{m_1}
          \\
          &
           T_2 \ar[ru]_{t_2} \ar[rd]_{l_2}
          &
          &
           {\txt{\scriptsize{$\theta_1 \!\!\Downarrow$}}}
          &
           W
          \\
          &
          &
           H_2 \ar[rd]_{m_2}
          \\
          &
          &
          &
           W_2 \ar[ruu]_{w_2}
         }
\xymatrix@R=11.5ex    
           {
            \\
            ,
            }
\hspace{3ex}
\xymatrix@R=1ex@C=4ex
         {
          &
          &
          &
           W_1 \ar[rdd]^{w_1}
          \\
          &
          &
           H_1 \ar[ru]^{n_1}
          \\
          &
           S_1 \ar[ru]^{r_1} \ar[rd]^{s_1}                    
          &
          &
           {\txt{\scriptsize{$\theta_2 \!\!\Downarrow$}}}
          &
           W 
          \\
           A \ar[ru]^{u_1} \ar[rd]_{u_2} 
          &
           {\txt{\scriptsize{$\beta \!\!\Downarrow$}}} 
          &
           S \ar[rd]^{n_2}
          \\
          &
           S_2 \ar[ru]_{s_2} \ar[rd]_{r_2}
          &
           {\txt{\scriptsize{$\phi_2 \!\!\Downarrow$}}}
          &
           W_2 \ar[ruu]_{w_2}
          \\
          &
          &
           H_2 \ar[ru]_{m_2}
         }
$$

%\pagebreak

The corresponding LL equation is :

\noindent
$
\xymatrix@R=1ex@C=4ex
         {
          &
           A \ar[rd]^{u_1}
          \\
          &
           {\txt{\scriptsize{$\xi_1 \!\!\Downarrow$}}}
          &
           S_1 \ar[rd]^{r_1}
          \\
           F \ar[r]^{y} \ar[rdd]_{z} \ar[ruu]^{x}
          &
           B \ar[ru]^{v_1} \ar[rd]_{h_1} 
          &
           {\txt{\scriptsize{$\gamma_1 \!\!\Downarrow$}}}
          &
           H_1 \ar[rd]^{n_1}
          \\
          &
           {\txt{\scriptsize{$\zeta_1 \!\!\Downarrow$}}}
          &
           T_1 \ar[ru]_{l_1} \ar[rd]_{m_1}
          &
           {\txt{\scriptsize{$\phi_1 \!\!\Downarrow$}}}
          &
           W_1 \ar[rdd]^{w_1}
          \\
          &
           C \ar[ru]_{k_1} \ar[rd]_{k_2}
          &
          {\txt{\scriptsize{$\varepsilon \!\!\Downarrow$}}}
          &
           T \ar[ru]_{m_1}
          \\
           &
           &
            T_2 \ar[ru]_{t_2} \ar[rd]_{l_2}
           &
           &
           {\txt{\scriptsize{$\theta_1 \!\!\Downarrow$}}}
           &
            W 
          \\
           &
           &
           &
            H_2\ar[rd]_{m_2}
          \\
           &
           &
           &
           &
            W_2 \ar[ruu]_{w_2}
          }
\hspace{6ex}
\xymatrix@R=18ex{\\ \txt{=}}
$

\nopagebreak

$
\hspace{14ex}
\xymatrix@R=18ex{\\ \txt{=}}
\hspace{6ex}
\xymatrix@R=1ex@C=4ex
         {
          &
          &
          &
          &
           W_1 \ar[rdd]^{w_1}
          \\
          &
          &
          &
           H_1 \ar[ru]^{n_1}
          \\
          &
          &
           S_1 \ar[ru]^{r_1} \ar[rd]^{s_1}                    
          &
          &
           {\txt{\scriptsize{$\theta_2 \!\!\Downarrow$}}}
          &
           W 
          \\
          &
           A \ar[ru]^{u_1} \ar[rd]^{u_2} 
          &
           {\txt{\scriptsize{$\beta \!\!\Downarrow$}}} 
          &
           S \ar[rd]^{n_2}
          \\
          &
           {\txt{\scriptsize{$\xi_2 \!\!\Downarrow$}}}
          &
           S_2 \ar[ru]^{s_2} \ar[rd]^{r_2}
          &
           {\txt{\scriptsize{$\phi_2 \!\!\Downarrow$}}}
          &
           W_2 \ar[ruu]_{w_2}
          \\
           F \ar[ruu]^{x} \ar[r]^{y} \ar[rdd]_{z}
          &
           B \ar[ru]^{v_2} \ar[rd]^{h_2}
          &
           {\txt{\scriptsize{$\gamma_2 \!\!\Downarrow$}}}
          &
           H_2 \ar[ru]_{m_2}
          \\
          &
           {\txt{\scriptsize{$\zeta_2 \!\!\Downarrow$}}}
          &
           T_2 \ar[ru]_{l_2}
          \\
          &
           C \ar[ru]_{k_2}
         }
$

\vspace{2ex}

\noindent To pass from the left side to the right side of this
equation use the LL equations of $(\varepsilon, \; \delta)$,
$(\theta_1, \; \theta_2)$ and $(\alpha, \; \beta)$, in this order. 
\end{proof}

\begin{definition}[equivalence of premorphisms] \label{equivalence}
Two premorphisms $\xi_1$,  $\xi_2$ are said to be \emph{equivalent}
when there exists a homotopy $(\alpha,\, \beta): \xi_1 \Rightarrow
\xi_2$. We shall write $\xi_1 \sim \xi_2$.
\end{definition}
Equivalence is indeed an equivalence relation. Proposition
 \ref{v.c.o.h} shows transitivity, the
 inverse  \mbox{2-cells} define a homotopy 
\mbox{$(\alpha^{-1}, \; \beta^{-1}): \xi_2 \Rightarrow \xi_1$} in the opposite
 direction, which shows
 symmetry, while reflexivity is obvious.

\vspace{1ex}

\pagebreak

\begin{definition}[Category $\cc{L}(F)$] \label{categoryLL} 
We define a category $\cc{L}(F)$ with objects pairs $(x, \; A)$, $x \in
FA$. Morphisms are equivalence classes of premorphisms, and composition
is defined by composing representative premorphisms. 
\end{definition}
It follows from
\ref{h.c.o.h} that composition is, up to equivalence, independent of the choice of
representatives, and
independent of the choice of the \mbox{2-cells} given by axiom F1 when
composing each pair of representatives. Since associativity holds and
identities exist, this construction actually does define a category.

\vspace{1ex}

%\end{document}

{\bf Some lemmas on pre-2-filtered categories}

\vspace{1ex} 

We establish now a few lemmas which are useful when proving the
fundamental properties of the construction LL.

\begin{lemma} \label{lemaF2}
Given a finite family of 2-cells 
$
\xymatrix@ur@R=2.5ex@C=2.5ex
         {
            E \ar[r]^f  \ar@{}[dr]|{\gamma_i \Downarrow}  \ar[d]_g 
          & A \ar[d]^{u_i}
          \\
            B  \ar[r]_{v_i} 
          & C_i 
         }
$, 
\mbox{$i = 1 \ldots n $,} 
there exists $A \mr{u} C$, $B \mr{v} C$, $C_i \mr{w_i} C, \; i = 1
\ldots n$, with invertible 
2-cells $\alpha_i , \; \beta_i$, such that the 2-cells
$$
\xymatrix@R=0.2ex@C=3ex
         {
          & A \ar[rdd]_{u_i} \ar@/^2ex/[rrrdd]^u
          \\
          && {\;\; \txt{\scriptsize{$\beta_i \!\Downarrow$}}} 
          \\
          E \ar[ruu]^f \ar[rdd]_g
          & {\txt{\scriptsize{$\gamma_i \!\Downarrow$}}}
          & C_i \ar[rr]^{w_i}
          && C
          \\
          && {\txt{\;\; \scriptsize{$\alpha_i \!\Downarrow$}}}
          \\
          & B \ar[ruu]^{v_i} \ar@/_2ex/[rrruu]_v 
         }
\hspace{2ex}
 \xymatrix@R=3.5ex
           {
            \\
            i = 1 \ldots n
            }     
$$
are all equal.

Given a second family of 2-cells  
$
\xymatrix@ur@R=2.5ex@C=2.5ex
         {
            H \ar[r]^h  \ar@{}[dr]|{\delta_i \Downarrow}  \ar[d]_l 
          & A \ar[d]^{u_i}
          \\
            B  \ar[r]_{v_i} 
          & C_i 
         }
$ (with same $u_i, \; v_i, \; C_i$),
we can assume that the same $u, \; v, \; w_i, \; \alpha_i , \;
\beta_i$, also equalize the 
2-cells of the second family.
\end{lemma}
\begin{proof}
Axiom F2 provides the case $n = 2$ with $u = w_1 u_1$, 
$v = w_2 v_2$, $\alpha_1 = \alpha$, $\beta_1 = id$, $\alpha_2 = id$, and 
$\beta_2 = \beta$. Using this case, induction is straightforward. For the second part,
if the 2-cells of the second family are not yet equalized, use the
lemma again (and patch the new 2-cells also into the first family)    
\end{proof}
From axiom F2 we deduce
\begin{lemma} \label{2-cellsallwaysequivalent}
Given any 2-cells 
$
\xymatrix@ur@R=2.5ex@C=2.5ex
         {
            E \ar[r]^f  \ar@{}[dr]|{\gamma_1 \Downarrow}  \ar[d]_g 
          & A \ar[d]^{u_1}
          \\
            B  \ar[r]_{v_1} 
          & C_1 
         }
 ,  \;
\xymatrix@ur@R=2.5ex@C=2.5ex
         {
            E \ar[r]^f  \ar@{}[dr]|{\gamma_2 \Downarrow}  \ar[d]_g 
          & A \ar[d]^{u_2}
          \\
            B  \ar[r]_{v_2} 
          & C_2 
         }
$
and an object $F \mr{x} E$, the premorphisms
$
\xymatrix@ur@R=2.5ex@C=2.5ex
         {
          F \ar[r]^{fx}  \ar@{}[dr]|{\!\!\!\gamma_1 x \Downarrow}  \ar[d]_{gx} 
          & A \ar[d]^{u_1}
          \\
            B  \ar[r]_{v_1} 
          & C_1 
         }
 ,  \;
\xymatrix@ur@R=2.5ex@C=2.5ex
         {
          F \ar[r]^{fx}  \ar@{}[dr]|{\!\!\!\gamma_2 x \Downarrow}  \ar[d]_{gx} 
          & A \ar[d]^{u_2}
          \\
            B  \ar[r]_{v_2} 
          & C_2 
         }
$
are equivalent.
\qed \end{lemma} 

From proposition \ref{basich} it follows that 
\begin{lemma} \label{basice}  
Given a pair of premorphism $\xi_1$, $\xi_2$ and a
pair of invertible two cells 
$\alpha$, $\beta$ that fit as follows, we have: 

\vspace{1ex}

$
\hspace{12ex}
\xymatrix@ur@R=2.5ex@C=2.5ex
         {
            F \ar[r]^x  \ar@{}[dr]|{\xi_1 \Downarrow}  \ar[d]_y 
          & A \ar[d]^{u_1}
          \\
            B  \ar[r]_{v_1} 
          & C_1
         }
\;\;\;\; \sim \;\;\;\;
\xymatrix@ur@R=2.5ex@C=2.5ex
         {
            F \ar[r]^x  \ar@{}[dr]|{\xi_1 \Downarrow}  \ar[d]_y 
          & A \ar[d]^{u_1}
          \\
            B \ar[r]_{\!\!{v_1}}  \ar@{}[dr]|{\alpha \Downarrow}
	                                                     \ar[d]_{{v_2}}  
          & C_1 \ar[d]^{w_1}
          \\
            C_2 \ar[r]_{{w_2}}
          & C
         }
$ \hspace{5ex} and 

$
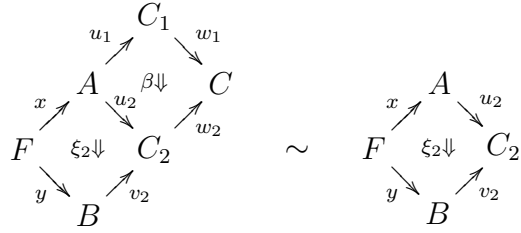
\hspace{12ex}
\xymatrix@ur@R=2.5ex@C=2.5ex
         {
            F \ar[r]^x  \ar@{}[dr]|{\xi_2 \Downarrow}  \ar[d]_y 
          & A \ar[r]^{u_1}  \ar@{}[dr]|{\beta \Downarrow}
	                                             \ar[d]^{\!\!\!\!{u_2}} 
          & C_1 \ar[d]^{w_1}
          \\
            B  \ar[r]_{{v_2}} 
          & C_2 \ar[r]_{{w_2}}
          & C
         }
\;\;\;\; \sim \;\;\;\;
\xymatrix@ur@R=2.5ex@C=2.5ex
         {
            F \ar[r]^x  \ar@{}[dr]|{\xi_2 \Downarrow}  \ar[d]_y 
          & A \ar[d]^{{u_2}}
          \\
            B  \ar[r]_{{v_2}} 
          & C_2
         }
$

\vspace{1ex}

When the pair $\alpha, \; \beta$ satisfy the LL equation, then
transitivity applied to these two equivalences yields the
equivalence $\xi_1 \sim \xi_2$.
\qed \end{lemma}
From this lemma and transitivity of equivalence we deduce
\begin{lemma} \label{ximaslejos}
Given a premorphism $\xi$, and a
pair of invertible two cells $\alpha$, $\beta$ that paste as follows, we have:
$$
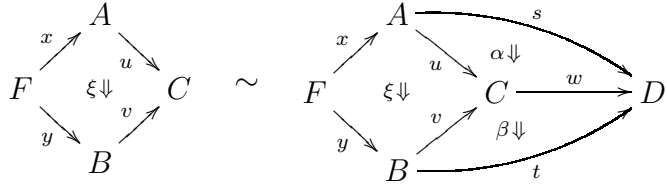
\xymatrix@R=0.75ex@C=2.75ex
         {
          & A \ar[rdd]_u 
          \\ 
          \\
          F \ar[ruu]^x \ar[rdd]_y
          & {\txt{\scriptsize{$\xi \!\Downarrow$}}}
          & C 
          \\
          \\
          & B \ar[ruu]^v  
         }           
\hspace{2ex}
\xymatrix@R=3.5ex
           {
            \\
            \sim
            }
\hspace{2ex}
\xymatrix@R=0.2ex@C=3ex
         {
          & A \ar[rdd]_u \ar@/^2ex/[rrrdd]^s
          \\
          && {\;\; \txt{\scriptsize{$\alpha \!\Downarrow$}}} 
          \\
          F \ar[ruu]^x \ar[rdd]_y
          & {\txt{\scriptsize{$\xi \!\Downarrow$}}}
          & C \ar[rr]^w
          && D
          \\
          && {\txt{\;\; \scriptsize{$\beta \!\Downarrow$}}}
          \\
          & B \ar[ruu]^v \ar@/_2ex/[rrruu]_t 
         }           
$$
\qed \end{lemma}     
   
\vspace{1ex}

{\bf The universal property of the construction LL}

\vspace{1ex}

A \emph{pseudocone}   for a 
2-functor $F$ with vertex the category $\cc{X}$ is a pseudonatural
transformation \mbox{$F \Mr{h} \cc{X}$} from $F$ to the 2-functor
which is constant at $\cc{X}$. Explicitly, it consists of a family of 
functors \mbox{$(h_A: FA \to \cc{X})_{A \, \in \cc{A}}$,} and a family
of invertible natural 
transformations \mbox{$(h_u: h_B \circ u \to h_A)_{(A \mr{u} B) \, \in
    \cc{A}}$.} 
 A morphism $h \Mr{\varphi} l$ of pseudocones (with same vertex) 
is a modification; as such, it consists of a
family of natural transformations 
$(h_A \Mr{\varphi_A} l_A)_{A \in \cc{A}}$. In accordance with notation
  \ref{abuse}, we have
$$
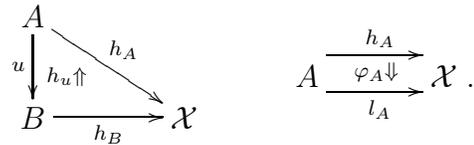
\xymatrix@C=8ex
          {
              A \ar[dr]^{h_A}  \ar[d]_{u}^(0.6){\; h_u \Uparrow} 
           \\
              B \ar[r]_{h_B} 
            & \cc{X}
           }
\hspace{6ex}
\xymatrix@C=7ex@R=2.4ex
         {
          \\
          A \ar@<1.5ex>[r]^{h_A} 
          \ar@{}@<-1ex>[r]^{\!\! \varphi_A \, \!\Downarrow}
	                                 \ar@<-1.2ex>[r]_{l_A} &
					 \cc{X}\; . 
         }
$$

This data is subject to the equations: 

\vspace{1ex}

\noindent PC0. $\hspace{8ex} h_{id_A} = id_{h_A}$

\noindent PC1.
$$
\xymatrix@C=8ex
          {
              A \ar[dr]^{h_A}  \ar[d]_{u}^(0.6){\; h_u \Uparrow} 
           \\
              B \ar[r]^{h_B} \ar[d]_{v}^(0.4){\; h_v \Uparrow}
            & \cc{X}
           \\ 
              C \ar[ur]_{h_C}
           }
\hspace{3ex} 
\xymatrix@R=6ex{\\ \txt{=}}
\hspace{3ex}
\xymatrix@C=8ex
          {
              A \ar[dr]^{h_A}  \ar[d]_{u} 
           \\
              B \ar@{}@<-0.5ex>[r]^(0.4){h_{vu} \Uparrow} \ar[d]_{v}
            & \cc{X}
           \\ 
              C \ar[ur]_{h_C}
           }
$$
PC2.
$$
\xymatrix@C=9ex
          {
              A \ar@<0.2ex>[dr]^{h_A}
                \ar@<-1.2ex>[d]^(0.4){\gamma}^(0.6){\! \Rightarrow}_{u}  
                \ar@<1.2ex>[d]^{\!v}^(0.6){\;\; h_v \Uparrow} 
           \\
              B \ar[r]_{h_B} 
            & \cc{X}
           }
\hspace{3ex} 
\xymatrix@R=3ex{\\ \txt{=}}
\hspace{3ex}
\xymatrix@C=8ex
          {
              A \ar[dr]^{h_A}  \ar[d]_{u}^(0.6){\; h_u \Uparrow} 
           \\
              B \ar[r]_{h_B} 
            & \cc{X}
           }
$$            
PCM.
$$
\xymatrix@C=8ex
          {
               A \ar@<-0.8ex>[dr]_(0.25){h_A\!\!\!\!\!}
                                ^(0.45){\!\!\varphi_A}^(0.6){\Uparrow}  
                \ar@<1.8ex>[dr]^{l_A} 
                \ar[d]_{u}^(0.8){\;\;\; h_u \Uparrow} 
           \\
              B \ar[r]_{h_B} 
            & \hspace{0ex} \cc{X}
           }
\hspace{3ex}
\xymatrix@R=3ex{\\ \txt{=}}
\hspace{3ex}
\xymatrix@C=8ex
          {
              A \ar@<0.8ex>[dr]^{l_A}  \ar[d]_{u}^(0.5){\; l_u \Uparrow} 
           \\
              B \ar@<1ex>[r]^{l_B} \ar@<-1.6ex>[r]_(0.5){h_B}^(0.5)
                                                {\varphi_B  \Uparrow}   
            & \cc{X}
           }
$$
Given a pseudo cone $F \Mr{h} \cc{Z}$ and a 2-functor 
$\cc{Z} \mr{s} \cc{X}$, it is clear and straightforward how to define
a pseudocone $F \Mr{sh} \cc{X}$ which is the composite 
$F \Mr{h} \cc{Z}  \mr{s} \cc{X}$; this determines a functor          
\mbox{$\cc{C}at[\cc{Z},\, \cc{X}] \times  \cc{P}\cc{C}[F,\, \cc{Z}] \mr{}
\cc{P}\cc{C}[F,\, \cc{X}]$.} 

\vspace{1ex}

\begin{theorem} \label{bicolimit}
Let $A \mr{u} B$ in $FA$ and $x \mr{\xi} y$ in $\cc{A}$. 
The following formulas define a pseudocone $F \Mr{\lambda} \cc{L}(F)$:
$$
\lambda_A(x) = F \mr{x} A \,, 
\hspace{3ex}
\lambda_A(\xi) =
\xymatrix@ur@R=2.5ex@C=2.5ex
         {
            F \ar[r]^x  \ar@{}[dr]|{\xi \Downarrow}  \ar[d]_y 
          & A \ar[d]^{{id}}
          \\
            A  \ar[r]_{{id}} 
          & A
         },
\hspace{3ex}
\lambda_u(x) = 
\xymatrix@ur@R=2.5ex@C=2.5ex
         {
            F \ar[r]^{ux}  \ar@{}[dr]|{id \Downarrow}  \ar[d]_x 
          & B \ar[d]^{{id}}
          \\
            A  \ar[r]_{{u}} 
          & B
         } 
$$
which induces by composition an equivalence of categories
$\cc{C}at[\cc{L}(F),\, \cc{X}] \mr{\cong} \cc{P}\cc{C}[F,\, \cc{X}]$, and
this equivalence is actually an \mbox{isomorphism.} We have:

\noindent for all $h$ there exists a unique $\stackrel{\sim}{h}$
such that $\stackrel{\sim}{h} \lambda = h$:
$
%For\; all\; h\; there\; exists\; \stackrel{\sim}{h}\;
%such\; that\; \; \stackrel{\sim}{h} \lambda = h\;:\;
\xymatrix
         {
            F \ar@2[r]^{\lambda} \ar@2[rd]_{\forall \, h}^{\; \equiv} 
          & \cc{L}(F) \ar@{{}{--}{>}}[d]^{\exists \, ! \, \stackrel{\sim}{h}}
          \\
          & \cc{X}
         }
$.
\end{theorem}
\begin{proof}
Functoriality of $\lambda_A$, naturality of $\lambda_u$ and equation
PC1 hold tautologically. The validity of PC2 means the following
equivalence of premorphisms:
$$
\xymatrix@R=1ex@C=4ex
         {
          &
           B \ar[rd]^{id}
          \\
          &
           {\txt{\scriptsize{$\gamma x \!\Downarrow$}}}
          &
           B \ar[rd]^{id}
          \\
           F \ar[r]^{vx} \ar[rdd]^x \ar[ruu]^{ux}
          &
           B \ar[ru]^{id} \ar[rd]^{id} 
          &
           {\txt{\scriptsize{$id  \!\Downarrow$}}}
          &
           B
          \\
          &
           {\txt{\scriptsize{$id \!\Downarrow$}}}
          &
           B \ar[ru]^{id}
          \\
          &
           A \ar[ru]^v
         }
\hspace{4ex}
\xymatrix@R=7ex{\\ \sim}
\hspace{0ex}
\xymatrix@ur@R=2.5ex@C=2.5ex
         {
          \\
          \\
            F \ar[r]^{ux}  \ar@{}[dr]|{id \Downarrow}  \ar[d]_x 
          & B \ar[d]^{{id}}
          \\
            A  \ar[r]_{{u}} 
          & B
         } 
$$
which is given by lemma \ref{2-cellsallwaysequivalent}.

We pass now to prove the universal property. Given $F \Mr{h} \cc{X}$,
define $\stackrel{\sim}{h}$ by the formulas: 
for
$  
\xymatrix@R=2.5ex@C=2.5ex@ur
         {
            F \ar[r]^x  \ar@{}[dr]|{\xi \Downarrow}  \ar[d]_y 
          & A \ar[d]^u
          \\
            B  \ar[r]_v 
          & C 
         }
$
 in $\cc{L}(F)$, 
$$
\xymatrix@R=3ex@C=3.5ex
           {
            \\
            {\stackrel{\sim}{h}(x) \;=\; (F} \ar[r]^(0.7){x}
            & A \ar[r]^(0.4){h_A}
            & \cc{X}),
            }
\hspace{1ex} 
\xymatrix@R=3ex
           {
            \\
           \stackrel{\sim}{h}(\xi) \; =
            }
\hspace{1ex}
\xymatrix@R=0.2ex@C=3ex
         {
          & A \ar[rdd]_u \ar@/^2ex/[rrrdd]^{h_A}
          \\
          && {\;\; \txt{\scriptsize{$h_{u}^{-1} \!\Downarrow$}}} 
          \\
          F \ar[ruu]^x \ar[rdd]_y
          & {\txt{\scriptsize{$\xi \!\Downarrow$}}}
          & C \ar[rr]^{h_C}
          && \cc{X}
          \\
          && {\txt{\;\; \scriptsize{$h_v \!\Downarrow$}}}
          \\
          & B \ar[ruu]^v \ar@/_2ex/[rrruu]_{h_B} 
         }
$$           
We have 
to show that the definition of $\stackrel{\sim}{h}(\xi)$ is compatible
with the equivalence of premorphisms. It is enough to consider the
two cases in lemma \ref{basice}. For the first case we have to show
the equation
$$
\xymatrix@R=2ex
           {
            \\
            (1)
            }
\hspace{9ex}
\xymatrix@R=3ex@C=3ex@ur
         {
            F \ar[r]^x  
              \ar@{}[dr]|{\xi \Downarrow}  
              \ar[d]_y 
          & A \ar[d]^{\!\!u} 
              \ar@/^6ex/[dd]^{\!\!h_A}
              \ar @/^3.5ex/@{{}{ }{}}[dd]|{h^{-1}_{u} \!\Downarrow}
          \\
            B  \ar[r]_v
               \ar@/_0.5ex/@{{}{ }{}}[dr]|(0.5){h_{v} \Downarrow}
               \ar@/_3ex/[dr]_(0.5){h_B}    
          & C  \ar[d]^{\!\!h_C}
          \\
          & \cc{X} 
         }
\hspace{2ex} 
\xymatrix@R=2ex
           {
            \\
             =
            }
\hspace{2ex}
\xymatrix@ur@R=3ex@C=3ex
         {
            F \ar[r]^x  
              \ar@{}[dr]|{\xi \Downarrow}  
              \ar[d]_y 
          & A \ar[d]^{\!\!u} 
              \ar@/^8ex/[ddd]^{\!\!h_A} 
              \ar @/^4ex/@{{}{ }{}}[ddd]|{h^{-1}_{wu} \Downarrow}
          \\
            B \ar[r]_{\!\!{v}}  
              \ar@{}[dr]|{\alpha \Downarrow}
              \ar[d]_{s} 
              \ar@/_4ex/@{{}{ }{}}[ddr]|(0.65){h_{rs} \Downarrow}
              \ar@/_8ex/[ddr]_(0.5){h_B}
          & C \ar[d]^{w}
          \\
            S \ar[r]_{\!\!r}
          & D \ar[d]^{\!\!\!h_D}
          \\
           & \cc{X}
         }
$$
Consider the following equation which follows from PC1:
$$
\xymatrix@R=2ex
           {
            \\
            (2)
            }
\hspace{5ex}
\xymatrix@ur@R=3ex@C=3ex
         {
            F \ar[r]^x  
              \ar@{}[dr]|{\xi \Downarrow}  
              \ar[d]_y 
          & A \ar[d]^{\!\!u} 
              \ar@/^8ex/[ddd]^{\!\!h_A} 
              \ar @/^5ex/@{{}{ }{}}[ddd]|(0.35){h^{-1}_{u} \!\Downarrow}
          \\
            B \ar[r]_{\!\!v}
              \ar@/_4.5ex/@{{}{ }{}}[ddr]|(0.5){h_{v} \Downarrow}
              \ar@/_8ex/[ddr]_(0.5){h_B}
          & C \ar[d]^{\!\!w} 
              \ar@/^5ex/[dd]^(0.4){\!\!h_C}
              \ar @/^3ex/@{{}{ }{}}[dd]|{\!\!h^{-1}_{w} \!\Downarrow}
              \ar @/_2ex/@{{}{ }{}}[dd]|(0.4){\!\!\!\!h_{w} \Downarrow}
              \ar@/_5ex/[dd]^(0.6){\!\!h_C}
          \\
          & D \ar[d]^{\!\!\!h_D}
          \\
           & \cc{X}
         }
\hspace{2ex} 
\xymatrix@R=2ex
           {
            \\
             =
            }
\hspace{2ex}
\xymatrix@ur@R=3ex@C=3ex
         {
            F \ar[r]^x  
              \ar@{}[dr]|{\xi \Downarrow}  
              \ar[d]_y 
          & A \ar[d]^{\!\!u} 
              \ar@/^8ex/[ddd]^{\!\!h_A} 
              \ar @/^4ex/@{{}{ }{}}[ddd]|{h^{-1}_{wu} \Downarrow}
          \\
            B \ar[r]_{\!\!v}
              \ar@/_4ex/@{{}{ }{}}[ddr]|(0.5){h_{wv} \Downarrow}
              \ar@/_8ex/[ddr]_(0.5){h_B}
          & C \ar[d]^{w}
          \\
          & D \ar[d]^{\!\!\!h_D}
          \\
           & \cc{X}
         }
$$
The right hand sides of equations (1) and (2)  are equal
by PC2, while the 
left hand sides are clearly equal. The second case in lemma
\ref{basice} is treated in a similar manner.
Functoriality of $\stackrel{\sim}{h}(\xi)$ follows from PC1 and PC2
with the same type of techniques as above. Finally, the equation
$\stackrel{\sim}{h} \lambda = h$ is immediate for the whole cone structure.       
\end{proof}

\section*{\normalsize 2. 2-Filtered 2-Categories} \label{fundamental}
\addtocounter{section}{1} \setcounter{equation}{0} 

%%%%%%%%%%%%%%%%%%%%%%%%%%%%%%% 

We refer here to Definition \ref{pre_2} of pre 2-filtered.

\begin{definition}
A 2-category $\cc{A}$ is defined to be \emph{pseudo 2-filtered} when  it satisfies the following three axioms:

\vspace{1ex}

\noindent FF1.  A stronger form of the axiom F1 of pre 2-filtered.
$$
Given \hspace{1ex}
\xymatrix@ur@R=2.5ex@C=2.5ex
         {
            E_1 \ar[r]^{f_1}   \ar[d]_{g_1} 
          & A
          \\
            B   
         } \;,
\hspace{0ex}
\xymatrix@ur@R=2.5ex@C=2.5ex
         {
            E_2 \ar[r]^{f_2}   \ar[d]_{g_2} 
          & A
          \\
            B   
         }
 \hspace{1ex} there \; exists \hspace{1ex}
\xymatrix@R=2.5ex@C=2.5ex@ur
         {
            E_1 \ar[r]^{f_1}  \ar@{}[dr]|{\gamma_1 \Downarrow}
	    \ar[d]_{g_1}  
          & A \ar[d]^u
          \\
            B  \ar[r]_v 
          & C 
         }
, \hspace{1ex}
\xymatrix@R=2.5ex@C=2.5ex@ur
         {
            E_2 \ar[r]^{f_2}  \ar@{}[dr]|{\gamma_2 \Downarrow}
	    \ar[d]_{g_2}  
          & A \ar[d]^u
          \\
            B  \ar[r]_v 
          & C 
         },
$$
$with \; \gamma_1 \; and \; \gamma_2 \; invertible \; 2 \text{-} \, cells$.

\vspace{1ex}

\noindent F2. Axiom F2 of pre 2-filtered.

\vspace{1ex}

\noindent F3.  Given two 2-cells as in axiom F2, with $B = A$, and  $u_1 = v_1$, $u_2 = v_2$, then there is a single $2$-cell 
$\varepsilon$ such that 
the LL-equation in F2 holds with $\varepsilon$ in place of both
$\alpha$ and $\beta$.
\end{definition}

\begin{remark} \label{BF2}
Given two $2$-cells as in axiom F2, with $B, \, C_1, \, C_2$ all equal to $A$, and $u_1,\, v_1,\, u_2,\,v_2$ all equal to $id_A$, then there exists $A \mr{w} C$ such that 
$w \,\gamma_1 = w \,\gamma_2$. Note that this Kennison axiom BF2, see Definition \ref{bifiltered}. 
\end{remark}
\begin{proof}
It follows immediately from axioms F2 and F3 that there is an invertible 
\mbox{$2$-cell} \mbox{$\varepsilon: w_1 \Mr{} w_2$,} such that 
$\varepsilon \,\gamma_1 = \varepsilon \,\gamma_2$. Cancelling 
$\varepsilon g$, we deduce that
\mbox{$w \,\gamma_1 = w \,\gamma_2$} with $w = w_1$.  
\end{proof}
 %%%%%%%%%%%%%%%%%%%%%%%%%%%
\begin{definition}
A 2-category $\cc{A}$ is defined to be \emph{2-filtered} when it is pseudo 2-filtered,
non empty, and satisfies in addition the following axiom.

\vspace{1ex}

\noindent F0.
$$
Given \hspace{5ex}
\xymatrix@ur@R=2ex@C=2ex
         {
          & A
          \\
            B   
         }
 \hspace{5ex} there \; exists \hspace{5ex}
\xymatrix@R=2.5ex@C=2.5ex@ur
         {
          & A \ar[d]^u
          \\
            B  \ar[r]_v 
          & C 
         }
\;.
$$
\end{definition}

When $\cc{A}$ is a trivial 2-category (the only 2-cells are the
identities), axiom F0 is the usual axiom in the definition of filtered
category, while our axiom FF1 is equivalent to the conjunction of
the two axioms PS1 and PS2 in the definition of pseudofiltered
category (cf \cite{G1} Expos\'{e} I). 

\vspace{1ex}

As was the case for axiom F1, in the presence of axiom WF3, axiom FF1 can
be replaced by the 
weaker version in which we do not require the 2-cells $\gamma_1$ and
$\gamma_2$ to be invertible.

\vspace{1ex} 

The following properties of the construction LL follow for pseudo
\mbox{2-filtered} 2-categories and not 
for pre 2-filtered 2-categories.
\begin{lemma} \label{u=v}
Given any morphism $(x,\,A) \to (y,\,A)$ in $\cc{L}(F)$, we can choose a
representative premorphism
$
\xymatrix@ur@R=2.5ex@C=2.5ex
         {
        F \ar[r]^{x}  \ar@{}[dr]|{\xi \Downarrow}  \ar[d]_{y} 
          & A \ar[d]^{u}
          \\
            A  \ar[r]_{v} 
          & C 
         }
$
with $u = v$. 
\end{lemma}
\begin{proof}
Consider 
$
\xymatrix@ur@R=2.5ex@C=2.5ex
         {
        F \ar[r]^{x}  \ar@{}[dr]|{\xi \Downarrow}  \ar[d]_{y} 
          & A \ar[d]^{u}
          \\
            A  \ar[r]_{v} 
          & C 
         }
$ 
and apply axiom FF1 to obtain invertible 2-cells $\alpha,\; \beta$
as follows:
$$
\xymatrix@R=0.75ex@C=2.75ex
         {
          & A \ar[rdd]^s 
          \\ 
          \\
          A \ar[ruu]^{id} \ar[rdd]_v
          & {\txt{\scriptsize{$\alpha \!\Downarrow$}}}
          & D 
          \\
          \\
          & C \ar[ruu]_t  
         }
\xymatrix@R=4.5ex
           {
            \\
            ,
            }
\hspace{1ex}
\xymatrix@R=0.75ex@C=2.75ex
         {
          & A \ar[rdd]^s 
          \\ 
          \\
          A \ar[ruu]^{id} \ar[rdd]_u
          & {\txt{\scriptsize{$\beta \!\Downarrow$}}}
          & D 
          \\
          \\
          & C \ar[ruu]_t  
         }
\xymatrix@R=3.75ex
           {
            \\
            :
            }
\hspace{4ex}
\xymatrix@R=0.2ex@C=3ex
         {
          & A \ar[rdd]_{u} \ar@/^2ex/[rrrdd]^s
          \\
          && {\;\; \txt{\scriptsize{$\beta \!\Downarrow$}}} 
          \\
          F \ar[ruu]^{x} \ar[rdd]_{y}
          & {\txt{\scriptsize{$\xi \!\Downarrow$}}}
          & C \ar[rr]^{t}
          && D
          \\
          && {\txt{\;\; \scriptsize{$\alpha^{-1} \!\Downarrow$}}}
          \\
          & A \ar[ruu]^{v} \ar@/_2ex/[rrruu]_s 
         }
$$ 
The proof follows by lemma \ref{ximaslejos}.
\end{proof}
\begin{lemma} \label{lemaF1}
Given a finite family of premorphisms 
$
\xymatrix@ur@R=2.5ex@C=2.5ex
         {
        F \ar[r]^{x_i}  \ar@{}[dr]|{\xi_i \Downarrow}  \ar[d]_{y_i} 
          & A \ar[d]^{u_i}
          \\
            B  \ar[r]_{v_i} 
          & C_i 
         }
$, 
\mbox{$i = 1 \ldots n $,} 
there exists $A \mr{u} C$, $B \mr{v} C$, $C_i \mr{w_i} C, \; i = 1
\ldots n$, with invertible 
2-cells $\alpha_i , \; \beta_i$ as in the diagram:
$$
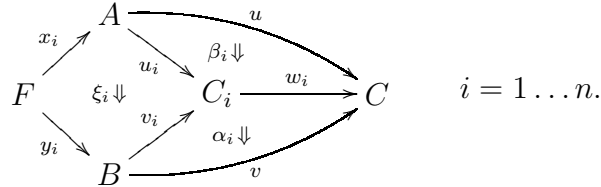
\xymatrix@R=0.2ex@C=3ex
         {
          & A \ar[rdd]_{u_i} \ar@/^2ex/[rrrdd]^u
          \\
          && {\;\; \txt{\scriptsize{$\beta_i \!\Downarrow$}}} 
          \\
          F \ar[ruu]^{x_i} \ar[rdd]_{y_i}
          & {\txt{\scriptsize{$\xi_i \!\Downarrow$}}}
          & C_i \ar[rr]^{w_i}
          && C
          \\
          && {\txt{\;\; \scriptsize{$\alpha_i \!\Downarrow$}}}
          \\
          & B \ar[ruu]^{v_i} \ar@/_2ex/[rrruu]_v 
         }
\hspace{4ex}
\xymatrix@R=3.5ex
           {
            \\
            i = 1 \ldots n.
            }
$$
When $A = B$ we can assume $u = v$.
\end{lemma}
\begin{proof}
Given $\xi_1,\; \xi_2$, apply FF1 to obtain invertible $\alpha$,
$\beta$ to fit as follows:
$$
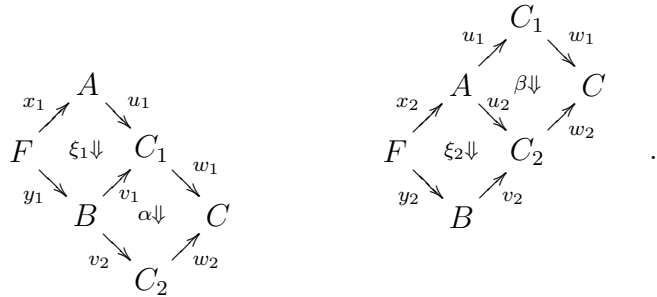
\xymatrix@ur@R=2.5ex@C=2.5ex
         {
          F \ar[r]^{x_1}  \ar@{}[dr]|{\xi_1 \Downarrow}  \ar[d]_{y_1} 
          & A \ar[d]^{u_1}
          \\
            B \ar[r]_{\!\!v_1}  \ar@{}[dr]|{\alpha \Downarrow}  \ar[d]_{v_2} 
          & C_1 \ar[d]^{w_1}
          \\
            C_2 \ar[r]_{w_2}
          & C
         }
\hspace{5ex}  \hspace{5ex}
\xymatrix@ur@R=2.5ex@C=2.5ex
         {
         F \ar[r]^{x_2}  \ar@{}[dr]|{\xi_2  \Downarrow}  \ar[d]_{y_2} 
          & A \ar[r]^{u_1}  \ar@{}[dr]|{\beta \Downarrow}
	    \ar[d]^{\!\!\!\!u_2} 
          & C_1 \ar[d]^{w_1}
          \\
            B  \ar[r]_{v_2} 
          & C_2 \ar[r]_{w_2}
          & C
         }
\;\;\;\;.
$$
This gives the case $n = 2$ with $u = w_1 u_1$, 
$v = w_2 v_2$, $\alpha_1 = \alpha$, $\beta_1 = id$, $\alpha_2 = id$, and 
$\beta_2 = \beta$. Using this case, induction is straightforward. For
the second part, do as in the proof of lemma \ref{u=v}. 
\end{proof}

%
%Notice that this lemma together with lemma \ref{ximaslejos} means that
%given any finite family of morphisms $(x,\;A) \to (y,\;B)$ in
%$\cc{L}(F)$, we can take representative premorphisms $(u,\;\xi_i,\;v)$
%all in the same home-set $\xi_i \in [ux,\;vy]$ of a single category $FC$.  
%

\begin{lemma} \label{fielmaslejos}
Given any pair of equivalent premorphisms
\mbox{$
\xymatrix@ur@R=2.5ex@C=2.5ex
         {
        F \ar[r]^{x}  \ar@{}[dr]|{\xi_1 \!\Downarrow}  \ar[d]_{y} 
          & A \ar[d]^{u_1}
          \\
            A  \ar[r]_{v_1} 
          & C_1 
         }
$
$\; \sim \;$
$
\xymatrix@ur@R=2.5ex@C=2.5ex
         {
        F \ar[r]^{x}  \ar@{}[dr]|{\xi_2 \!\Downarrow}  \ar[d]_{y} 
          & A \ar[d]^{u_2}
          \\
            A  \ar[r]_{v_2} 
          & C_2 
         }
$}
, if $u_1 = v_1$ and $u_2 = v_2$, then we can choose a homotopy defined by a single (invertible) $2$-cell $\varepsilon$, 
$w_1 \Mr{\varepsilon} w_2$, 
 $(\varepsilon,\, \varepsilon): \xi_1 \Rightarrow \xi_2$.
\begin{proof}
It follows immediately from axioms F2 and F3.
\end{proof}
\end{lemma}

\begin{lemma} \label{fielmaslejos2}
Given two arrows
$\xymatrix{x \ar@<1.25ex>[r]^{\xi_1} \ar@<-1.25ex>[r]^{\xi_2} & y}$
 in $FA$, if 
$\lambda_A(\xi_1) = \lambda_A(\xi_2)$ in $\cc{L}(F)$, then there exists $A \mr{w} C$ such that $w \,\xi_1 = w \,\xi_2$ in $FC$. 
\end{lemma}
%%%%%%%%%%%%%%%%%%%%%%%%%
\begin{proof}
Recall the definition of $\lambda_A$ in Theorem \ref{bicolimit}. By the previous Lemma with $C_1 = C_2=C$, and $u_1,\, v_1,\, u_2,\,v_2$ all equal to $id_C$, it follows that there is a $2$-cell $w_1 \Mr{\varepsilon} w_2$, such that 
$\varepsilon \,\xi_1 = \varepsilon \,\xi_2$. Since $\varepsilon$ is invertible, it follows that $w_1 \,\xi_1 = w_1 \,\xi_2$. Compare with Remark \ref{BF2}.
\end{proof}
%%%%%%%%%%%%%%%%%%%%%%%%% 
 
\begin{theorem} \label{finitepresentation}
Let $\cc{A}$ be a pre 2-filtered 2-category, 
\mbox{$F: \cc{A} \mr{} \cc{C}at$} a \mbox{2-functor,} and $\cc{P}$ a finite category.
Consider the 2-functor 
\mbox{$F^\cc{P}: \cc{A} \mr{} \cc{C}at$} defined by $F^\cc{P}(A) =
(FA)^\cc{P}$, and the canonical functor:
$$
\Diamond: \cc{L}(F^\cc{P}) \mr{} \cc{L}(F)^\cc{P}
\hspace{4ex} (given \; by \; theorem \; \ref{bicolimit}).
$$
Then, $\Diamond$  is an equivalence of categories provided that
{$\cc{A}$ is 2-filtered} or that 
{$\cc{A}$ is pseudo 2-filtered and $\cc{P}$ is connected}. 
\end{theorem}
\begin{proof}
Notice that an object $F^\cc{P} \mr{} A$ in $\cc{L}(F^\cc{P})$ is by
definition a
diagram $\cc{P} \mr{} FA$. We shall prove, in turn, that $\Diamond$ is (a)
essentially surjective, 
(b) faithful, and (c) full. 

\vspace{1ex}

\noindent (a) \emph{essentially surjective}: 
We shall see that given a diagram \mbox{$\cc{P} \to \cc{L}(F)$,} there
exists $A \in \cc{A}$ and a factorization (up to isomorphism):
$$
\xymatrix
        {
         & \cc{P} \ar[d] \ar[ld]^{\;\;\cong}
         \\
         FA \ar[r]^{\lambda_A}
         & \cc{L}(F)
        } 
$$
Consider explicitly an object in
$\cc{L}(F)^\cc{P}$:
$$
F \mr{x_k} A_k ,\; k \in \cc{P}, \hspace{2ex} 
\xymatrix@R=2.5ex@C=2.5ex@ur
         {
            F \ar[r]^{x_p}  \ar@{}[dr]|{\varphi_f \Downarrow}
	    \ar[d]_{x_q}  
          & A_p \ar[d]^{u_f}
          \\
            A_q  \ar[r]_{v_f} 
          & A_f 
         },
\hspace{2ex} p \mr{f} q  \in  \cc{P}. 
$$
satisfying equations $ \varphi_{f \circ g} \,\sim\, \varphi_f \circ_{\gamma} \varphi_g$ 
for all composable pairs $f,\; g$. 

Let $\cc{Q} \subset \cc{P}$ be a part of $\cc{P}$ for which   
there exists $A$, $w_k$, $\psi_f$ such that:
$$
F \mr{x_k} A_k ,\; k \in \cc{Q}, \hspace{1ex}
\xymatrix@R=2.5ex@C=2.5ex@ur
         {
            F \ar[r]^{x_p}  \ar@{}[dr]|{\varphi_f \Downarrow}
	    \ar[d]_{x_q}  
          & A_p \ar[d]^{u_f}
          \\
            A_q  \ar[r]_{v_f} 
          & A_f 
         }
\hspace{1ex} \sim \hspace{1ex}
\xymatrix@R=2.5ex@C=2.5ex@ur
         {
            F \ar[r]^{x_p}  \ar@{}[dr]|{\psi_f \Downarrow}
	    \ar[d]_{x_q}  
          & A_p \ar[d]^{w_p}
          \\
            A_q  \ar[r]_{w_q} 
          & A 
         },
\hspace{2ex} p \mr{f} q  \in  \cc{Q}. 
$$
$\cc{Q}$ is not necessarily a subcategory, but we agree that if \mbox{$p
\mr{f} q \, \in \cc{Q}$,} then we consider $p \in \cc{Q}$ and $q  \in
\cc{Q}$.
The equations  $ \psi_{f \circ g} \,\sim\, \psi_f \circ_{\gamma}
\psi_g$ hold for all composable pairs $f,\; g$ in $\cc{Q}$ with $f
\circ g$ also in 
$\cc{Q}$. By lemma \ref{fielmaslejos2} we can assume strict
equality $ \psi_{f \circ g} \,=\, \psi_f \circ \psi_g$ in $FA$. 

We shall see that if $p \mr{g} q$ is not in $\cc{Q}$, we
can add it to $\cc{Q}$ in such a way that the enlarged part retains
the same property. In what it follows we use without mention lemma
\ref{ximaslejos}.   

1) $p \in \cc{Q}, \; q \notin \cc{Q}$ or $q \in \cc{Q}, \; p \notin
   \cc{Q}$: In the first case apply axiom F1 to obtain an 
   invertible 2-cell $\beta$ as follows:
$$
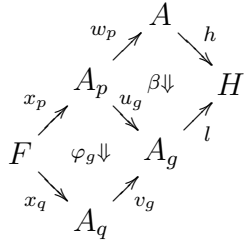
\xymatrix@ur@R=2.5ex@C=2.5ex
         {
            F \ar[r]^{x_p}  \ar@{}[dr]|{\varphi_g  \Downarrow}  \ar[d]_{x_q} 
          & A_p \ar[r]^{w_p}  \ar@{}[dr]|{\beta \Downarrow}
	    \ar[d]^{\!\!\!\!u_g} 
          & A \ar[d]^{h}
          \\
            A_q  \ar[r]_{v_g} 
          & A_g \ar[r]_{l}
          & H
         }
$$
The new $A$ is
$H$, the new $w_k$ are $hw_k$, all $k \in \cc{Q}$, the 
new $\psi_f$ are $h\psi_f$, all $f \in \cc{Q}$, and, finally, 
$w_q = lv_g$, and $\psi_g$ is the 2-cell above. The other case can be
proved in the same way. 

2) $p \in \cc{Q}, \; q \in \cc{Q}$: Apply axiom
   FF1 to obtain 
   invertible 2-cells $\alpha,\; \beta$ as follows:
$$
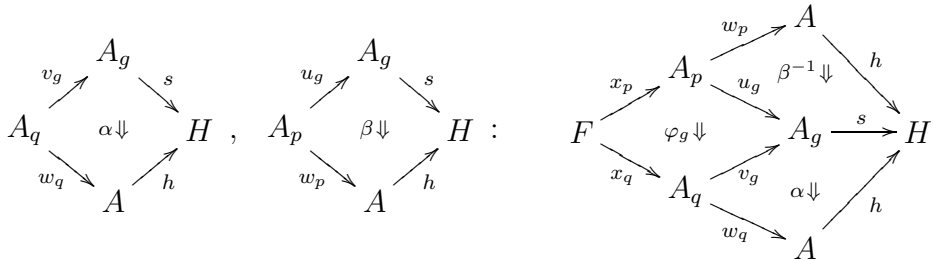
\xymatrix@R=0.5ex@C=2.75ex
         {
          \\
          & A_g \ar[rdd]^s 
          \\ 
          \\
          A_q \ar[ruu]^{v_g} \ar[rdd]_{w_q}
          & {\txt{\scriptsize{$\alpha \!\Downarrow$}}}
          & H 
          \\
          \\
          & A \ar[ruu]_h  
         }
\xymatrix@R=7.5ex
           {
            \\
            ,
            }
\hspace{1ex}
\xymatrix@R=0.5ex@C=2.75ex
         {
          \\
          & A_g \ar[rdd]^s 
          \\ 
          \\
          A_p \ar[ruu]^{u_g} \ar[rdd]_{w_p}
          & {\txt{\scriptsize{$\beta \!\Downarrow$}}}
          & H 
          \\
          \\
          & A \ar[ruu]_h  
         }
\xymatrix@R=6.75ex
           {
            \\
            :
            }
\hspace{4ex} 
\xymatrix@R=1ex@C=4ex
         {
          && A \ar[ddr]^h
          \\
          & A_p \ar[ur]^{w_p} \ar[dr]^{\!\!u_g}
          & {\txt{\scriptsize{$\beta^{-1} \!\Downarrow$}}}
          \\
          F \ar[ur]^{x_p} \ar[dr]_{x_q}
          & {\txt{\scriptsize{$\varphi_{g} \!\Downarrow$}}}
          & A_g \ar[r]^s
          & H
          \\
          & A_q \ar[ur]_{\!\!v_g} \ar[dr]_{w_q}
          &{\txt{\scriptsize{$\alpha \!\Downarrow$}}}
          \\
          && A \ar[uur]_h
         }
$$
The new $A$ is $H$, the new $w_k$ are $hw_k$, all $k \in \cc{Q}$, the
new $\psi_f$ are $h\psi_f$, all $f \in \cc{Q}$, and, finally, $\psi_g$
is the 2-cell above on the right. This proof holds whether $p \neq  q$ or
$p = q$.

3) $p \notin \cc{Q}, \; q \notin \cc{Q}$:
Apply axiom F0 to obtain $A_g \mr{w} H$, $A \mr{h} H$. The new $A$ is
$H$, the new $w_k$ are $hw_k$, all $k \in \cc{Q}$, the 
new $\psi_f$ are $h\psi_f$, all $f \in \cc{Q}$, and, finally, 
$w_p = lu_g$, $w_q = lv_g$, $\psi_g = l \varphi_g$. If $p = q$, use lemma
\ref{u=v} to assume $u_g = v_g$, and in this way $w_p$ is uniquely defined. 

Notice that if $\cc{P}$ is connected, it is not necessary to consider
this case since we can always choose $p \mr{g} q$ such that either
$p$, or $q$, or both, are in $\cc{Q}$. Thus for connected $\cc{P}$
axiom F0 is not necessary.

It is clear that any singleton $\{(k,\; f = id_k)\}$ serves as an initial
$\cc{Q}$, thus we can assume $\cc{Q} = \cc{P}$. 

To finish the proof observe that given any arrow $ux \mr{\psi} vy$ in $FA$,
the square 
$
\xymatrix@C=10ex
        {
          (ux,\:A) \ar[r]^{(id,\,id,\,u)}_\cong \ar[d]^{\lambda_{A}(\psi)}
        & (x,\:A)  \ar[d]^{(u,\,\psi,\,v)}
        \\
          (vy,\:A) \ar[r]^{(id,\,id,\,v)}_\cong
        & (y,\:B)
        }
$
\mbox{commutes} in $\cc{L}(F)$.

Notice that if $\cc{P}$ is a poset, case 2) cannot happen, so axiom
FF1 is not necessary. For posets $\cc{P}$ the functor $\Diamond$ is
essentially surjective also for pre 2-filtered 2-categories.
\vspace{1ex}

\noindent (b) \emph{faithful:} Consider two premorphisms 
$
\xymatrix@ur@R=2.5ex@C=2.5ex
         {
        F^\cc{P} \ar[r]^{x}  \ar@{}[dr]|{\xi \Downarrow}  \ar[d]_{y} 
          & A \ar[d]^{u}
          \\
            B  \ar[r]_{v} 
          & C 
         }
, \;
\xymatrix@ur@R=2.5ex@C=2.5ex
         {
        F^\cc{P} \ar[r]^{x}  \ar@{}[dr]|{\eta \Downarrow}  \ar[d]_{y} 
          & A \ar[d]^{s}
          \\
            B  \ar[r]_{t} 
          & D 
         }
$
in $\cc{L}(F^\cc{P})$,  
$
\xymatrix@ur@R=2.5ex@C=2.5ex
         {
        F \ar[r]^{x_k}  \ar@{}[dr]|{\xi_k \Downarrow}  \ar[d]_{y_k} 
          & A \ar[d]^{u}
          \\
            B  \ar[r]_{v} 
          & C 
         }
, \;
\xymatrix@ur@R=2.5ex@C=2.5ex
         {
        F \ar[r]^{x_k}  \ar@{}[dr]|{\eta_k \Downarrow}  \ar[d]_{y_k} 
          & A \ar[d]^{s}
          \\
            B  \ar[r]_{t} 
          & D 
         }
$, $k \in \cc{P}$.
To be equivalent in $\cc{L}(F)^\cc{P}$ means that there are homotopies
$(\alpha_k,\, \beta_k) : \xi_k \Rightarrow \eta_k$ given by invertible
2-cells  
$
\xymatrix@ur@R=2.5ex@C=2.5ex
         {
        B \ar[r]^{v}  \ar@{}[dr]|{\alpha_k \Downarrow}  \ar[d]_{t} 
          & C \ar[d]^{w_k}
          \\
            D  \ar[r]_{h_k} 
          & H_k 
         }
, \;
\xymatrix@ur@R=2.5ex@C=2.5ex
         {
        A \ar[r]^{u}  \ar@{}[dr]|{\beta_k \Downarrow}  \ar[d]_{s} 
          & C \ar[d]^{w_k}
          \\
            D  \ar[r]_{h_k} 
          & H_k 
         }
$, $k \in \cc{P}$.
From lemma \ref{lemaF2} it readily follows that we can assume there are
single invertible 2-cells $\alpha$, $\beta$ which define all the
homotopies  $(\alpha,\, \beta) : \xi_k \Rightarrow \eta_k$. But this
means that $\xi$ and $\eta$ are equivalent in $\cc{L}(F^\cc{P})$.

Notice that this proof of the faithfulness of the functor $\Diamond$
holds for pre 2-filtered 2-categories. 

\vspace{1ex}

\noindent (c) \emph{full}: Consider two objects 
$F^\cc{P} \mr{x} A$, $F^\cc{P} \mr{y} B$  in $\cc{L}(F^\cc{P})$. 
\mbox{A premorphism} in $\cc{L}(F)^\cc{P}$ consists of a family
$
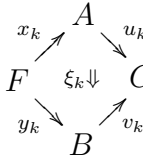
\xymatrix@R=2.5ex@C=2.5ex@ur
         {
            F \ar[r]^{x_k}  \ar@{}[dr]|{\xi_k \Downarrow}
	    \ar[d]_{y_k}  
          & A \ar[d]^{u_k}
          \\
            B  \ar[r]_{v_k} 
          & C_k 
         },
\; k \in \cc{P}. 
$
From lemma  \ref{lemaF1} we can assume all the $u_k,\, v_k, \, C_k$
to be equal to a single $u,\, v, \, C$. But this is the data for a
premorphism in $\cc{L}(F^\cc{P})$. For the naturality equations we
proceed as in the proof of \mbox{faithfulness in (b).}

Here axiom FF1 is inevitable (lemma \ref{lemaF1}), and plays the role
of axiom PS2 in the filtered category case. Notice that a function
between sets viewed as a functor between trivial categories is
injective precisely when it is a full functor. 
\end{proof}

\vspace{1ex}

We state now an important corollary of theorem
\ref{finitepresentation}. Let $\cc{C}at_{f\ell}$ be the 2-category
of finitely complete categories and finite limit preserving
functors. We have: 
\begin{theorem}
Let $\cc{A}$ be a 2-filtered 2-category, and $\cc{A} \mr{F}
\cc{C}at_{f\ell}$ a 2-functor. Then, the category $\cc{L}(F)$ has
finite limits, the pseudocone functors $FA \mr{\lambda_A} \cc{L}(F)$
preserve finite limits and induce an equivalence of categories
$\cc{C}at_{f\ell}[\cc{L}(F),\, \cc{X}] \mr{\cong}
\cc{P}\cc{C}_{f\ell}[F,\, \cc{X}]$; this equivalence is actually an
isomorphism. 
\qed \end{theorem} 

\vspace{1ex}

{\bf Kennison notion of bifiltered 2-category}

\vspace{1ex}

In \cite{K}, Kennison considers the following notion:

\begin{definition}[Kennison] \label{bifiltered}
A 2-category $\cc{A}$ is defined to be \emph{bifiltered} when it
satisfies the following three axioms:

\vspace{1ex}

\noindent BF0.
Given two objects $A$, $B$, there exists $C$ and $A \to C$, $B \to C$.

\vspace{1ex}

\noindent BF1.
Given two arrows  
$\xymatrix{A \ar@<1.25ex>[r]^{f} \ar@<-1.25ex>[r]^{g} & B}$, there 
  exists $B \mr{u} C$ 
and  an  invertible  2-cell $\gamma: \, uf \cong ug$.

\vspace{1ex}

\noindent BF2. Given two 2-cells
$\xymatrix@C=10ex{A \ar@<1.5ex>[r] 
\ar@<-1.5ex>[r]^{\gamma_1 \!\Downarrow \, \gamma_2 \!\Downarrow} &
B}$, 
there exists $B \mr{u} C$ such that $u\gamma_1 = u\gamma_2$.
 \end{definition}

This notion of bifiltered 2-category is equivalent to our notion of
\mbox{2-filtered 2-category.} The proof of this is elementary although not
entirely trivial, and we leave it as an interesting exercise.

\vfill\eject

{\bf Eduardo J. Dubuc }

\vspace{1ex}

Departamento de Matem\'atica,

F. C. E. y N., Universidad de Buenos Aires,

1428 Buenos Aires, Argentina.

\vspace{2ex}

{\bf Ross Street}

\vspace{1ex}

Department of Mathematics,

Macquarie University,

Sydney, Australia.

\end{document}